\documentclass{article}
\usepackage{amsmath,amssymb,amsthm,graphicx}

\usepackage[dvips]{color}

\theoremstyle{plain}
\newtheorem{thm}{Theorem}
\newtheorem{lem}[thm]{Lemma}
\newtheorem{prop}[thm]{Proposition}
\newtheorem{cor}[thm]{Corollary}

\theoremstyle{definition}

\newtheorem*{qtn}{Question}

\theoremstyle{remark}
\newtheorem*{rem}{Remark}
\newtheorem*{claim}{Claim}

\newcommand{\N}{\mathrm{N}}
\newcommand{\D}{\mathrm{D}}
\newcommand{\E}{\mathrm{E}}

\title{Strongly invertible knots, rational-fold branched coverings and hyperbolic spatial graphs}

\author{Kazuhiro {\sc Ichihara} \vspace{\medskipamount}\\
School of Mathematics Education\\
Nara University of Education\\
Takabatake-cho\\
Nara 630-8528, JAPAN\\
{\tt ichihara@nara-edu.ac.jp}
\vspace{\bigskipamount}
\and Akira {\sc Ushijima} \vspace{\medskipamount}\\
Division of Mathematical and Physical Sciences\\
Graduate School of Natural Science and Technology\\
Kanazawa University\\
Kanazawa 920-1192, JAPAN\\
{\tt ushijima@kenroku.kanazawa-u.ac.jp}
}

\begin{document} \maketitle
\begin{abstract}
A construction of a spatial graph from a strongly invertible knot 
was developed by the second author, and
a necessary and sufficient condition for the given spatial graph to be hyperbolic
was provided as well. 
The condition is improved in this paper.
This enable us to show that certain classes of knots can yield hyperbolic spatial graphs via the construction.

\vspace{\medskipamount}

\begin{flushleft}

{\bf Key words:} spatial graph, theta-curve, strongly invertible knot, simple knot, hyperbolic manifold.\\

\vspace{\smallskipamount}

{\bf 2000 Mathematics Subject Classifications:}  
Primary: 57M50; secondary: 05C10.

\end{flushleft}
\vspace{\bigskipamount}
\end{abstract}
%
%
%
%
%
%
%
%
%
%
%
%
\section{Introduction}
\label{sec_intro}

Throughout this paper the category we will work with is differentiable.
We start this section with recalling several basic terminologies frequently used in this paper.
Let $M$ be a compact orientable three-dimensional manifold.
If it has boundary, then we denote the boundary by $\partial M$.
A {\em spatial graph}\/ in $M$ means 
the image of an embedding of a graph,
which will be regarded as a one-dimensional cellular complex as usual.
We often identify a spatial graph and the graph of its source in this paper;
a vertex, an edge, etc.\ of a spatial graph $G$ 
refers that of the graph appearing as the source of $G$.
Spatial graphs are assumed to be finite throughout this paper.

A spatial graph is especially called a {\em link}\/
if it is homeomorphic to a disjoint union of circles,
and a link is particularly called a {\em knot}\/ if its connected component is one.
Let $\N (G)$ be a regular neighborhood of a spatial graph $G$ in $M$,
and let $\E (G)$ be the {\em exterior}\/ of $G$ in $M$,
namely $\E (G)$ is the closure of $M - \N (G)$ in $M$.
Suppose that any connected components of $\partial M$ are 
closed orientable surfaces except spheres.
Then a spatial graph $G$ in $M$ is said to be {\em hyperbolic\/}
if $\E (G)$ minus all the toral boundary components
admits a complete hyperbolic structure of finite volume with geodesic boundary.

Let $G$ be a spatial graph in the three-dimensional sphere $S^3$,
and $\varphi$ an orientation-preserving self-diffeomorphism of the pair $(S^3 , G)$.
For a positive integer $n$,
we denote by $\varphi^n$ the map obtained by iterating $\varphi$ $n$ times.
Then $\varphi$ is said to be {\em periodic of order $n$}
when  $\varphi^m$ is not the identity map on $(S^3 , G)$
for any integer $m$ with $1 \leq m \leq n-1$,
and $\varphi^n$ is the identity map.
It is known that if such $\varphi$ has a fixed point,
then the fixed point set is homeomorphic to a circle.
See Chapter~1 of \cite{mhb84} for example.
A spatial graph $G$ is said to be {\em strongly periodic of order $n$}\/
if there is a periodic self-diffeomorphism $\varphi$ of order $n$
such that the fixed point set of $\varphi$
intersects with each connected component of $G$ at exactly two points,
and that $n$ is the maximal period that $(S^3 ,G)$ can admit with the same fixed point set.
Such a self-diffeomorphism is called a {\em strong inversion}\/ when $n=2$,
and strongly periodic links of order $2$ are especially called 
{\em strongly invertible}\/ links.

In \cite{us_HypGraph}, the second author has developed a method, called
{\em the $n / 2$-fold cyclic branched covering}\/, 
to construct strongly periodic spatial graphs in $S^3$ 
from a strongly invertible link $L$ in $S^3$.
Let $\alpha$ be the axis of a strong inversion $\iota$ of the pair $(S^3,L)$.
Then the construction consists of the following two steps:
\begin{description}
\item[Step~1.]
	Take a quotient of $(S^3,L)$ by $\iota$.
	Then we have $S^3 / \iota$, which is homeomorphic to $S^3$.
	In $S^3 / \iota$, we have
	a spatial graph $G$ as the quotient of $L$ by $\iota$,
	which has two univalent vertices on $\alpha / \iota = \alpha$.
\item[Step~2.]
	Take the $n$-fold cyclic branched covering of $(S^3 / \iota, G)$ along $\alpha$.
	Then we have the $n$-fold cyclic branched covering of $S^3 / \iota$,
	which is again homeomorphic to $S^3$.
	In the $n$-fold cyclic branched covering of $S^3 / \iota$,
	we have a spatial graph $G'$ as the pre-image of $G$ by the covering.
\end{description}
%
%
%
%
%
%
%
%
%
%
\begin{figure}[ht]
        \begin{center}
		\includegraphics[width=100mm]{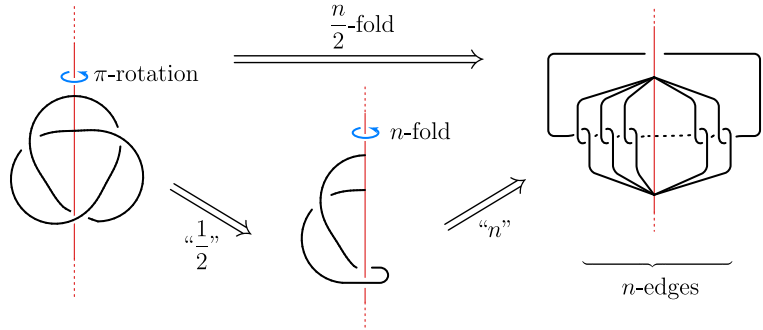}
        \end{center}
        \caption{The process of constructing ${\frac{n}{2}}$-fold cyclic branched covering of the trefoil knot}
        \label{f_n_over_2}
\end{figure}

We call the spatial graph $G'$ obtained by this procedure 
the {\em $n / 2$-fold cyclic branched covering of $L$ along $\alpha$}.
A spatial graph is called a {\em $\theta_n$-curve}\/ for some integer $n \geq 3$
when each connected component has exactly two vertices of valency $n$.
If we start with a strongly invertible link,
then the $n / 2$-fold cyclic branched covering
yields strongly periodic $\theta_n$-curves.
We have arcs as the $1/2$-fold cyclic branched covering when $n=1$,
and we then call the union of them as a strongly periodic $\theta_1$-curve in this paper.
Similarly we have the original $L$ as the $2/2$-fold cyclic branched covering when $n=2$,
and we occasionally call it a strongly periodic $\theta_2$-curve.

%
%
%
%
%
%
%
%
%
%
%
\section{When can we have strongly periodic hyperbolic $\theta_n$-curves?}

Let $M$ be a three-dimensional manifold
admitting a periodic self-diffeomorphism $\varphi$ of order $n$.
Then a submanifold $S$ of $M$ 
is said to be {\em equivariant under $\varphi$}
if either $\varphi^i (S) = S$ or $\varphi^i (S) \cap S = \emptyset$
holds for any integer $i$ with $1 \leq i \leq n-1$.
Let $\Gamma$ be a finite group acting on $M$.
Then $S$ is said to be {\em $\Gamma$-equivariant}\/
if it is equivariant under the action of any element of $\Gamma$.

Let $K$ be a non-trivial strongly invertible knot in ${S}^3$ 
admitting a strong inversion with axis $\alpha$, and 
$\theta_n^K$ the strongly periodic $\theta_n$-curve 
obtained from $K$ by the $n / 2$-fold cyclic branched covering along $\alpha$.
Then it has been essentially proved in \cite{us_HypGraph} that 
$\theta_n^K$ is hyperbolic for any $n \geq 3$
if and only if 
$\E ( K \cup \alpha )$ contains no essential tori
which are equivariant under the action of the strong inversion.

%
%
%
%
%
%
%
%
%
\subsection{Main theorem}

The first purpose of this paper is to simplify and improve the result given in \cite{us_HypGraph}.
A spatial graph in $S^3$ is said to be {\em trivial}\/
when it can be ambient isotopic into an embedded two-dimensional sphere.
%
%
%
%
%
%
\begin{thm} \label{thm_equivalent}
Let $K$ be a non-trivial strongly invertible knot in $S^3$
admitting a strong inversion with axis $\alpha$, and 
$\theta_n^K$ the strongly periodic $\theta_n$-curve
obtained from $K$ by the $n / 2$-fold cyclic branched covering along $\alpha$.
Then the following conditions are equivalent:
\begin{enumerate}
\item \label{enu_all}
	For any $n \geq 3$, $\theta_n^K$ is hyperbolic.
\item \label{enu_some}
	For some $n \geq 3$, $\theta_n^K$ is hyperbolic.
\item \label{enu_atoroidal}
	Any embedded torus in $\E (K \cup \alpha)$ is compressible in $\E (K \cup \alpha)$.
\item \label{enu_equatoroidal}
	Any embedded equivariant torus in $\E (K \cup \alpha)$ is compressible in $\E (K \cup \alpha)$.
\item  \label{enu_atolong}
	Any embedded torus in $\E (K)$ which is disjoint from $\alpha$ is compressible in $\E (K)$.
\item \label{enu_equatolong}
	Any embedded equivariant torus in $\E (K)$ which is disjoint from $\alpha$ is compressible in $\E (K)$.
\end{enumerate}
\end{thm}

\vspace{\bigskipamount}

\begin{figure}[ht]
        \begin{center}
		\includegraphics[width=30mm]{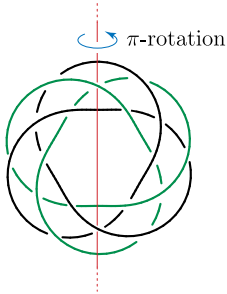}
        \end{center}
        \caption{Two-component torus link and its axis of a strong inversion}
        \label{f_T46.eps}
\end{figure}

There are some remarks about this theorem.
\begin{itemize}
\item
	We shall use standard language in three-dimensional topology.
	We here include some of them for reader's convenience. 
	See \cite{jab80} on such terminologies for example. 
	Let $M$ be a compact three-dimensional manifold.
	
	A properly embedded surface $F$ in $M$ other than a sphere or a disc 
	is said to be {\em compressible}\/ in $M$
	if it admits a disc $D$ such that $F \cap D = F \cap \partial D = \partial D$ 
	and that $\partial D$ is not null-homotopic on $F$. 
	Such a disc is called  a {\em compressing disc}\/ for $F$. 
	Otherwise $F$ is said to be  {\em incompressible}. 
	
	Analogously, a disc $D'$ embedded in $M$ 
	is called a {\em boundary-compressing disk}\/ for $F$ 
	if $\partial D' = a \cup b$,
	where $a := F \cap D' = F \cap \partial D'$
	is an arc which is not parallel to $\partial F$ in $F$
	and $b := \partial M \cap D' = \partial M \cap \partial D'$
	is an arc in $\partial M$ which is not homotopic into $\partial F$ relative $\partial b$.
	If there is no boundary-compressing disk for $F$, 
	then $F$ is said to be {\em boundary-incompressible}. 
	
	The manifold $M$ is said to be {\em irreducible}\/ 
	if any embedded sphere bounds an embedded ball,
	and {\em boundary-irreducible}\/
	if each connected component of $\partial M$ is incompressible.
	
	Two embedded surfaces $F$ and $F'$ are said to be {\em parallel}\/
	if they bound a product region in $M$. 
	In particular, when $F$ is parallel to a subsurface of $\partial M$, 
	it is said to be {\em boundary-parallel}.
	
	Finally, an embedded surface is said to be {\em essential}\/
	if it is incompressible, boundary-incompressible, and not boundary-parallel.  
\item
	Though we do not mention $\alpha$ in the notation $\theta_n^K$,
	the spatial graph $\theta_n^K$ depends not only on both $K$ and $n$ 
	but also on the choice of the axis $\alpha$.
	See Subsection~\ref{subsec_TwoAxes} for such an example.
\item
	Since $\partial \E ( K \cup \alpha )$
	consists of an orientable closed surface of genus three,
	there is no embedded torus which is boundary-parallel. 
	Thus Condition~\eqref{enu_atoroidal} means that $\E ( K \cup \alpha )$ is {\em atoroidal}.
\item
	The assumption that $K$ is a {\em knot}\/ is essential.
	Suppose that $K$ is a two-component non-trivial torus link
	and $\alpha$ is the axis of the strong inversion.
	See Figure~\ref{f_T46.eps} for example.
	Then $K$ is known to be strongly invertible,
	and $\E ( K \cup \alpha )$ is known to be atoroidal.
	But it has been proved as Theorem~1.3 in \cite{us_HypGraph}
	that $\theta_n^K$ is not hyperbolic.
\end{itemize}

%
%
%
%
%
%
%
%
%
%

%
%
%
%
%
%
%
%
%
%
%
%
\subsection{Equivariant Torus Theorem} \label{subsec_ett}
In this subsection we give our version of the so-called ``Equivariant Torus Theorem". 
It should be emphasized that
the boundary of the manifold in the proposition below is not assumed to be incompressible.

\begin{prop} \label{prop_ett}
For a compact orientable irreducible three-dimensional manifold with non-empty boundary,
if it admits an action of a finite group and 
if it contains an essential torus,
then it contains an embedded equivariant incompressible torus as well.
\end{prop}

To prove this proposition, we prepare a definition. 
Let $M$ be an orientable three-dimensional manifold, 
and let $S$ be an essential torus 
properly embedded in $M$.
Then we say that $S$ is said to be {\em canonical}\/
if any other properly embedded essential torus or annulus can be isotoped to be disjoint from $S$.
%
%
%
%
%
%
%
%
%
%
%
\begin{lem} \label{lem_canonical}
For a compact orientable irreducible three-dimensional manifold,
if it admits an action of a finite group and if it contains a canonical torus,
then the torus is isotoped to be equivariant under the action.
\end{lem}
\begin{proof}
Let $M$ be a compact orientable irreducible three-dimensional manifold
admitting an action of a finite group $\Gamma$, and $T$ a canonical torus in $M$.
Then we can choose a $\Gamma$-invariant Riemannian metric on $M$
whose boundary has non-negative mean curvature.
Such a metric can be constructed by averaging any metric over the action.
Then we have a least area torus $T'$ isotopic to $T$ by virtue of the result in \cite{fhs83}.
Precisely we use an extension of Theorems~1.1 and 1.2 of \cite{fhs83}
to the non-empty boundary case, 
which is explained at the 3rd paragraph in page 635 of \cite{fhs83}.
Since the torus $T'$ is also canonical as is $T$,
its images by the action are isotopically disjoint or coincident, 
and thus they are actually distinct or coincident by 
an extension of Theorem 6.2 of \cite{fhs83} to the non-empty boundary case, 
which is also explained at the 3rd paragraph in page 635 of \cite{fhs83}. 
This indicates that $T'$ is equivariant under the action. 
\end{proof}

%
%
%
%
%
%
%
%
\begin{proof}[Proof of Proposition~\ref{prop_ett}]
Let $M$ be a compact orientable irreducible three-dimensional manifold
admitting an action of a finite group, and $T$ an essential torus in $M$.

If $M$ contains a canonical torus,
then it is isotoped to be equivariant under the action by Lemma~\ref{lem_canonical}.
So we assume that $M$ has no canonical torus from now on.

Let us set $M' := M$ if $\partial M$ is incompressible.
If, on the other hand, $\partial M$ is compressible, we set $M'$ as follows.
We consider a maximal set of compression discs for $\partial M$, and
compress $\partial M$ along the discs.
Then we have a sub-manifold $N$, possibly disconnected, in $M$.
By standard cut-and-paste argument together with the assumption of the irreducibility for $M$,
any incompressible surface in $M$ can be isotoped apart from $\partial N$. So
we can find a torus $T'$ in $N$ which is isotopic to $T$ in $M$.
Then we set $M'$ as the connected component of $N$ containing $T'$.

Here we note that any toral component of $\partial M'$ is boundary-parallel in $M$.
For, as mentioned above,
any incompressible surface in $M$ can be isotoped
apart from $\partial M' \subset \partial N$.
So the toral component must be boundary-parallel in $M$
by the assumptions that there is no canonical torus in $M$ 
and that $\partial M'$ is incompressible in $M$ by its construction.
This implies that $T'$ is essential in $M'$
since $T'$ is not boundary-parallel in $M$ as is $T$.

As a result we thus find an irreducible three-dimensional submanifold $M'$
containing an essential torus $T'$ but not containing any canonical torus.
Then $M'$ must be either a Seifert fibred space or an $I$-bundle
by Proposition~3.2 of \cite{nsp97}.
Thus all connected components of $\partial M'$ are tori
and all of them are, as claimed above, boundary-parallel in $M$.
This concludes that $M'=M$, 
and thus $M$ contains an incompressible torus as a boundary component,
which gives an embedded equivariant incompressible torus.
This completes the proof of Proposition~\ref{prop_ett}.
\end{proof}

%
%
%
%
%
%
%
%
%
%
\subsection{Lemmata for the proof of Theorem~\ref{thm_equivalent}}
\label{subsec_preparation}

In this subsection, we will prepare three lemmata to prove Theorem~\ref{thm_equivalent}. 

The first one would be well-known in three-dimensional manifold theory.
Actually a similar result is shown as Theorem in \cite{glb84}.
%
%
%
%
%
%
%
%
%
%
\begin{lem} \label{lem_ECD}
Suppose that an irreducible three-dimensional manifold
admitting an action of a finite cyclic group 
contains an embedded torus equivariant under the action.
If the torus is compressible,
then there exists a compression disc which is equivariant under the action.
\end{lem}

\begin{proof}
Let $M$ be an irreducible three-dimensional manifold
admitting a non-trivial action $\varphi$ of period $n$,
and let $T$ be
an embedded torus in $M$ which is equivariant under $\varphi$
with compression disc $D$.

We first suppose that
$D$ has non-empty intersection with another $\varphi^i ( T )$
for some $i \in \{ 1, 2, \dotsc , n-1 \}$.
Since $\varphi^i (T)$ are mutually disjoint,
the intersection appearing on $D$ must consist of circles in the interior of $D$.
Then, since $M$ is irreducible,
one can reduce the circles which are inessential on the torus by isotopy.
Take an inner-most one on $D$,
and replace $T$ and $D$
with $\varphi^i ( T )$ and the chosen disc.
Thus, without loss of generality,
we can assume that 
there is a compression disc for such a torus
having no intersection with any other tori $\varphi^i (T)$,
and we denote a pair of the torus and the compression disc
by $T$ and $D$ again.

Take an equivariant regular neighborhood $\N(T)$ of $T$ under $\varphi$,
and we define a submanifold $W$ in $M$ as follows:
\begin{equation*}
W := M - \bigcup_{i=1}^n \varphi^i ( \N(T) ) .
\end{equation*}
Here we take $\N(T)$ sufficiently small 
so that $D \cap W$ gives a compression disc in $W$ for $\partial_0 W$,
where $\partial_0 W$ is a connected component of $\partial W$
given by one of two connected components of $\partial \N (T)$.
We then apply the Equivariant Loop Theorem
(see \cite{my84ic1})
to $\left( W, \partial W, D \cap W \right)$
so that we can find a new compression disc $D'$ for $\partial_0 W$ in $W$,
which is equivariant under $\varphi$. 
Enlarging this disc and we can find a compression disc for $T$ in $M$,
which is equivariant under $\varphi$.
\end{proof}

We can prove the following two lemmata using this lemma.

%
%
%
%
%
%
%
%
%
\begin{lem} \label{lem_covering}
Let $L$ be a strongly invertible link in $S^3$
admitting a strong inversion with axis $\alpha$,
and $\theta_n^L$ the spatial graph with axis $\alpha_n$
which is obtained from $L$ by $n / 2$-fold cyclic branched covering along $\alpha$.
Then, if any embedded equivariant torus in $\E( \theta_n^L \cup \alpha_n )$
is compressible in $\E( \theta_n^L \cup \alpha_n )$
for some $n \geq 1$,
then so is it
for all $n \geq 1$.
\end{lem}
\begin{rem}
Strongly periodic action on $\E( \theta_1^L \cup \alpha_1 )$ means the trivial action.
So an embedded {\em equivariant}\/ torus in $\E( \theta_1^L \cup \alpha_1 )$
is nothing but an embedded torus,
or in other words,
any embedded torus is an equivariant torus in $\E( \theta_1^L \cup \alpha_1 )$ in our sense.
\end{rem}
%
%
%
%
%
%
\begin{proof}
Our proof consists of the following two steps:
\begin{description}
\item[Step 1.]
	Suppose that any embedded equivariant torus in $\E(\theta_n^L \cup \alpha_n )$
	is compressible for some $n \geq 2$,
	and prove that so is it in $\E(\theta_1^L \cup \alpha_1 )$.
\item[Step 2.]
	Suppose that any embedded equivariant torus in $\E(\theta_1^L \cup \alpha_1 )$
	is compressible,
	and prove that so is it in $\E(\theta_n^L \cup \alpha_n )$ for any $n \geq 2$.
\end{description}
By taking a regular neighborhood suitably, 
we can assume that $\E(\theta_n^L \cup \alpha_n)$ is 
an $n$-fold cyclic covering space of $\E(\theta_1^L \cup \alpha_1)$.

To prove the first step,
let $T$ be an embedded (equivariant) torus in $\E(\theta_1^L \cup \alpha_1)$,
and we will find its compression disc in $\E(\theta_1^L \cup \alpha_1)$.
Consider the preimage of $T$ in $\E(\theta_n^L \cup \alpha_n)$ 
by the covering projection.
It consists of either an equivariant torus, or a tuple of mutually disjoint equivariant tori.
Let $\widetilde{T}$ be one of its connected components.
Then the assumption says that
there is a compression disc for $\widetilde{T}$ in $\E(\theta_n^L \cup \alpha_n)$.
Furthermore $\E(\theta_n^L \cup \alpha_n)$ is irreducible
since so is $S^3$ and $\E(\theta_n^L \cup \alpha_n) \subset S^3$.
So we can apply Lemma~\ref{lem_ECD}
and find an equivariant compression disc for $\widetilde{T}$ in $\E(\theta_n^L \cup \alpha_n)$,
and the projection of the disc in $\E(\theta_1^L \cup \alpha_1)$
becomes the compression disc for $T$ in $\E(\theta_1^L \cup \alpha_1)$.
We have thus proved the first step.

To prove the second step,
let $T$ be an embedded equivariant torus in $\E(\theta_n^L  \cup \alpha_n)$,
and we will find its compression disc in $\E(\theta_n^L  \cup \alpha_n)$.
Since $T$ has no intersection with the axis $\alpha_n$ of $\varphi$, 
the projection of $T$ by the quotient map induced by the action of $\varphi$
is an embedded torus in $\E(\theta_1^L \cup \alpha_1)$.
Since we are assuming that
any embedded torus in $\E(\theta_1^L \cup \alpha_1)$ is compressible,
there is a compression disc for the projection of $T$.
Since the disc has no intersection with the branched set $\alpha_1$,
its lift in $\E(\theta_n^L  \cup \alpha_n)$ is a compression disc for $T$.
We have thus proved the second step,
which means that the proof of Lemma~\ref{lem_covering} have been completed.
\end{proof}

%
%
%
%
%
%
%
%
%
%
\begin{lem} \label{lem_comp}
Let $G$ be a strongly periodic spatial graph in $S^3$
admitting a periodic action with axis $\alpha$.
If an embedded equivariant torus in $\E (G)$
which is disjoint from $\alpha$ is compressible in $\E (G)$,
then it has a compression disc which is disjoint from $\alpha$ as well.
\end{lem}
%
%
%
%
%
%
%
%
%
%
\begin{proof}
Let $\varphi$ be the periodic of order $n$ on $(S^3 , G)$ with axis $\alpha$.
By taking a regular neighborhood of $G$ suitably,
we can assume that $\E (G)$ admits the action of $\varphi$.
Let $T$ be an embedded torus in $\E (G)$
which is equivariant under $\varphi$. 
Now we have its compression disc in $\E (G)$ by the assumption of this lemma.
Since $\E (G)$ is irreducible as is $S^3$,
we can find an equivariant compression disc $D$ for $T$ by Lemma~\ref{lem_ECD}.

If $D$ has no intersection with $\alpha$, 
then we are done.
So we assume $D \cap \alpha \ne \emptyset$ from now on.
Since $\alpha$ is the fixed point set of the action on $(S^3 , G)$
and since $D$ is equivariant under this action,
we can say that $\varphi^i (D) = D$ holds for any $i \in \{ 1, 2 , \dotsc , n \}$,
i.e., $\left. \varphi \right|_D$ can be regarded as a periodic diffeomorphism on $D$.
Then,
by the well-known classical result, 
there is an element $r$ in the orthogonal group of dimension two
together with a homeomorphism $h$ on $D$
such that $\left. \varphi \right|_D = h \, r \, h^{-1}$.
See Theorem~3.1 in \cite{ck94} for its modern proof.
So we have the following two cases,
depending on the orientation of $\left. \varphi \right|_D$ on $D$:
\begin{itemize}
\item $\left. \varphi \right|_D$, and thus $r$, is orientation-reversing, or
\item $\left. \varphi \right|_D$, and thus $r$, is orientation-preserving.
\end{itemize}

If the first case happens,
the fixed point set of $r$ is a diameter of $D$ 
and thus the fixed point set of $\left. \varphi \right|_D$
is an embedded arc on $D$,
which is $D \cap \alpha$.
We take a sufficiently small regular neighborhood $\N (\alpha)$ in $\E (G)$
so that $D \cap \partial \N (\alpha)$ consists of two arcs parallel to $D \cap \alpha$.
Then replace $D \cap \N (\alpha)$ by 
one of two simply-connected regions on $\partial \N (\alpha)$
bounded by four arcs:
two of them come from $T \cap (\partial \N (\alpha) - D)$
and the other two arcs are parallel ones appearing as $D \cap \partial \N (\alpha)$.
We thus have a compression disc for $T$
which has no intersection with $\alpha$.

We lastly consider the second case.
Since $r$ is non-trivial,
the set of fixed points for $r$ consists of a single point,
and so does the set of fixed points for $\left. \varphi \right|_D$.
This means that $D \cap \alpha$ consists of a single point.
The axis $\alpha$ is the unknot by the affirmative solution of the Smith Conjecture.
See page 4 of \cite{mhb84}.
So we can regard $T$ as the boundary of a trivially-embedded solid torus $V$
with core curve $\alpha$.
Since $G \cup \alpha$ is connected by the definition of strongly periodic $\theta_n$-curve,
and since $T \cap \alpha = \emptyset$,
$V$ entirely contains $G \cup \alpha$.
So the closure of $\E ( G \cup \alpha ) - V$ is also a trivially-embedded solid torus.
Therefore we can find a compression disc for $T = \partial V$ 
in $\E( G \cup \alpha )-V \subset \E(G \cup \alpha)$,
which induces a compression disc for $T$ in $\E(G \cup \alpha)$.
\end{proof}

\subsection{Proof of Theorem~\ref{thm_equivalent}}
Now we are in the position to give a proof of Theorem~\ref{thm_equivalent}.
%
%
%
%
%
%
%
%
%
%
%
\begin{proof}[Proof of Theorem~\ref{thm_equivalent}]
Though the equivalence between \eqref{enu_all} and \eqref{enu_equatoroidal}
has been partially proved as Theorem~1.2 in \cite{us_HypGraph},
we will give a whole proof of Theorem~\ref{thm_equivalent} for completeness.
We show it following the diagram below:
\begin{equation*}
\begin{matrix}
\eqref{enu_all} &&&& \eqref{enu_atoroidal}\\
& \nwarrow && \nearrow && \searrow\\
\downarrow && \eqref{enu_equatoroidal} &&&& \eqref{enu_atolong}\\
& \nearrow && \nwarrow && \swarrow\\
\eqref{enu_some} &&&& \eqref{enu_equatolong}
\end{matrix}.
\end{equation*}
\begin{description}
\item[$\eqref{enu_all} \rightarrow \eqref{enu_some}.$]
	This is clear.
\item[$\eqref{enu_some} \rightarrow \eqref{enu_equatoroidal}.$]
	Assume that	the strongly periodic $\theta_n$-curve $\theta_n^K$
	is hyperbolic for some $n \geq 3$.
	Then there is no essential torus in $\E (\theta_n^K)$
	by Thurston's study of geometric structure on three-dimensional manifolds.
	So any embedded equivariant torus in $\E (\theta_n^K)$ which is disjoint from $\alpha_n$
	is compressible in $\E (\theta_n^K)$ as well.
	Then
	there is a compression disc in $\E (\theta_n^K)$ which has no intersection with $\alpha_n$
	by Lemma~\ref{lem_comp}.
	This means that any embedded equivariant torus in $\E(\theta_n^K \cup \alpha_n)$
	is compressible in $\E(\theta_n^K \cup \alpha_n)$.
	So we can apply Lemma~\ref{lem_covering} for $n=2$
	and we thus have \eqref{enu_equatoroidal}
	since $\E(\theta_2^K \cup \alpha_2) = \E(K \cup \alpha)$.
\item[$\eqref{enu_equatoroidal} \rightarrow \eqref{enu_all}.$]
	Since $\theta_1^K$ is an arc whose terminal
	vertices lie on $\alpha_1$,
	the manifold $\E (\theta_1^K)$ is topologically the three-dimensional ball,
	and $\alpha_1$ turns out to be two arcs in the ball,
	which we will denote by $\alpha_1$ again.
	Thus $( \E (\theta_1^K) , \alpha_1 )$ can be regarded
	as a two-string tangle.
	This realization is the so-called Montesinos trick.
	See \cite{moi75, bha96}.
	
	Since $K$ is a non-trivial knot in $S^3$, 
	the exterior of $K$ in $S^3$ is irreducible and boundary-irreducible.
	The two-string tangle $( \E (\theta_1^K) , \alpha_1 )$
	is thus prime by Theorem~3.5.17 in \cite{ka96b}.
	
	Moreover
	Lemma~\ref{lem_covering} can be applied by the assumption
	so that
	any embedded equivariant torus in $\E(\theta_1^K \cup \alpha_1)$
	is compressible in $\E(\theta_1^K \cup \alpha_1)$.
	Since the finite group acting on $\E(\theta_1^K \cup \alpha_1)$ is the trivial group,
	the result above actually means that
	any embedded torus in $\E(\theta_1^K \cup \alpha_1)$
	is compressible in $\E(\theta_1^K \cup \alpha_1)$.
	This implies that the tangle $( \E (\theta_1^K) , \alpha_1 )$ atoroidal,
	or it is the so-called {\em hyperbolic}\/ tangle in the sense of \cite{ka96b}.

	Let $( \D_{\E (\theta_1^K)} , \D_{\alpha_1} )$
	be the double of $( \E (\theta_1^K) , \alpha_1 )$
	along the sphere $\partial \E (\theta_1^K)$.
	Then $\D_{\alpha_1}$ gives a hyperbolic link in $S^3$
	by Exercise~3.6.4 together with Theorem~3.6.6 both in \cite{ka96b}.
	Now, as we saw in Theorem~2.1 in \cite{us_HypGraph},
	the $n$-fold cyclic branched covering of 
	$( \D_{\E (\theta_1^K)} , \D_{\alpha_1} )$ along $\D_{\alpha_1}$
	is a closed hyperbolic manifold for any $n \geq 3$.
	
	Each manifold obtained by the covering admits an involution 
	along the surface in the manifold induced from $\partial \E (\theta_1^K)$.
	The action can be regarded as an isometry
	by Mostow-Prasad's rigidity theorem (see \cite{fra04} for example).
	Thus the quotient by this action is a hyperbolic manifold with totally geodesic boundary.
	As a result, the $n$-fold cyclic branched covering of $( \E (\theta_1^K) ,\alpha_1 )$
	along $\alpha_1$
	is a hyperbolic manifold with totally geodesic boundary for any $n \geq 3$.
	This conclusion means that \eqref{enu_all} holds
	by the definition of $n/2$-fold cyclic branched covering.
\item[$\eqref{enu_atoroidal} \rightarrow \eqref{enu_atolong}.$]
	This is trivial since $\E (K \cup \alpha)$ is a submanifold in $\E (K)$.
\item[$\eqref{enu_atolong} \rightarrow \eqref{enu_equatolong}.$]
	This is trivial.
\item[$\eqref{enu_equatolong} \rightarrow \eqref{enu_equatoroidal}.$]
	This is an immediate consequence of Lemma~\ref{lem_comp}.
\item[$\eqref{enu_equatoroidal} \rightarrow \eqref{enu_atoroidal}.$]
	Proposition~\ref{prop_ett} guarantees the contraposition of this statement as follows:
	suppose that there is an incompressible torus in $\E (K \cup \alpha)$.
	Since the genus of $\partial \E (K \cup \alpha)$ is three,
	it is not boundary-parallel, i.e., it is actually essential.
	Then
	we can find an embedded equivariant incompressible torus
	in $\E (K \cup \alpha)$ by Proposition~\ref{prop_ett}.
\end{description}
We have thus completed the proof of Theorem~\ref{thm_equivalent}.
\end{proof}

Let $G_n$ be a hyperbolic strongly periodic $\theta_n$-curve in $S^3$.
Then, as in the previous theorem, 
a knot can be obtained from $G_n$
by the $2/n$-fold cyclic branched covering;
namely first take the quotient by the periodic diffeomorphism of order $n$,
and then take the double branched covering along the axis of the periodic diffeomorphism.
The knot obtained by this procedure
must be a non-trivial strongly invertible knot,
and the spatial graph obtained from this knot by $l/2$-fold cyclic branched covering
is a hyperbolic strongly periodic $\theta_l$-curve,
which can be called the {\em spatial periodic $\theta_l$-curve
obtained from $G_n$ by the $l/n$-fold cyclic branched covering}.
We have thus obtained the following corollary
as an immediate consequence of Theorem~\ref{thm_equivalent}:
%
%
%
%
%
%
\begin{cor} \label{cor_hypgraph}
Let $G_n$ be a hyperbolic strongly periodic $\theta_n$-curve in $S^3$.
Then, for any integer $l \geq 3$, 
$G_l$ is a hyperbolic strongly periodic $\theta_l$-curve in $S^3$ as well,
where $G_l$ is 
constructed from $G_n$
by the $l/n$-fold cyclic branched covering. \qed
\end{cor}

Now, once we have obtained Theorem~\ref{thm_equivalent},
the following question naturally arises:
\begin{qtn}
Which strongly invertible knot satisfies the conditions of Theorem~\ref{thm_equivalent}?
\end{qtn}
We will give several partial answers to this question in the next section.

%
%
%
%
%
%
%
%
%
%
%
%
%
%
%
%
\section{Which strongly invertible knot can be a source of strongly periodic hyperbolic $\theta_n$-curves?} \label{sec_ex}

We start this section with proving the following lemma.
%
%
%
%
%
%
%
%
%
\begin{lem} \label{lem_parallel}
Let $G$ be a strongly periodic spatial graph in $S^3$
admitting a periodic action with axis $\alpha$.
Suppose that there is an embedded equivariant torus in $\E (G)$
which is disjoint from $\alpha$
and boundary-parallel.
Then $G$ must be a knot, and the torus has a compression disc in $\E (G)$;
this implies that $G$ is the trivial knot.
\end{lem}
%
%
%
%
%
%
%
%
%
\begin{proof}
The assumption of this lemma says that
there is an embedded torus, say $T$, in $\E (G \cup \alpha)$ 
which is parallel to a boundary component 
$\partial_0 \E (G)$ of $\partial \E (G)$ in $\E (G)$.
This means that $G$ has at least one connected component
which is homeomorphic to a circle.
Thus, by the definition of strongly periodic spatial graph,
the strongly periodic action must be a strong inversion,
which we denote by $\iota$.
We now start this proof with showing that $G$ is actually a knot.

By taking regular neighborhoods of $G$ and $\alpha$ suitably,
we can assume that 
both $\E (G)$ and $\E (G \cup \alpha)$ 
admit an action induced by the periodic diffeomorphism.
Since $T$ is parallel to $\partial_0 \E (G)$,
there is a compact submanifold $V$ in $\E (G)$ such 
that $\partial V$ consists of $T$ and $\partial_0 \E (G)$
and that $V$
is homeomorphic to
the product manifold
$T^2 \times [0,1]$ constructed by the two-dimensional torus $T^2$.
Since $T$ can be regarded as in $S^3$, 
it separates $S^3$ into two compact submanifolds.
We denote one of them containing $V$ by $M_V$.
Then, since $\partial_0 \E (G) \cap \alpha \ne \emptyset$ and $T \cap \alpha = \emptyset$
whereas $G \cup \alpha$ is connected, 
$M_V$ entirely contains $G \cup \alpha$.
Thus all boundary components of $\E (G)$ must be contained in $M_V$.
However, 
since $M_V$ consist of the union of a solid torus and $V$,
$M_V$ can contain at most one component of $G$. 
This implies that $G$ must be a knot,
and thus $\partial_0 \E (G) = \partial \E (G)$.

We next show that the torus $T$ has a compression disc in $\E (G)$.
Since $T$ is equivariant under $\iota$,
possible arrangements between $T$ and $\iota (T)$ are
either $\iota (T) = T$ or $\iota (T) \cap T = \emptyset$.
So one of the following two cases can happen for $\iota (T)$:
\begin{itemize}
\item $\iota (T) \subset V$, or
\item $\iota (T) \subset \E (G) - V$.
\end{itemize}

We first consider the latter case.
Then, since $\iota$ is a homeomorphism which preserves $\partial \E (G)$,
the submanifold $\iota (V)$ is again homeomorphic to $T^2 \times [0,1]$ with
$\partial \iota (V) = \iota (T) \cup \partial \E (G)$.
Thus $\iota (T)$ is again parallel to $\partial \E (G)$.
Furthermore, since $\iota = \iota^{-1}$ and $\iota (T) \cap T = \emptyset$ in this case,
the embedded torus $T$ must be contained in $\iota (V)$.
So we can regard $\iota (V)$ and $\iota (T)$ as $V$ and $T$,
and the situation turns out to be the former case.

We then consider the former case.
Recall that what we want to prove is the existence of a compression disc for $T$ in $\E (G)$,
and we will do it by contradiction.
So we assume that $T$ is incompressible in $\E (G)$.
Then, since $\iota$ is a homeomorphism,
$\iota (T)$ is also incompressible in $\E (G)$,
and thus especially incompressible in $V$.
Since $V$ has a homeomorphism $\psi$ to $T^2 \times [0,1]$,
the embedded torus $\psi \circ \iota (T)$ is also incompressible in $T^2 \times [0,1]$.
Then, by LEMMA~(1) in \cite{hei81}, any incompressible surface in $T^2 \times [0,1]$
must be parallel to both $T^2 \times \{ 0 \}$ and $T^2 \times \{ 1 \}$,
and thus so is $\psi \circ \iota (T)$.
Since $\partial \E (G)$ is preserved by $\iota$,
we have $\psi \circ \iota (V) \subseteq T^2 \times [0,1]$,
which is equivalent to $\iota (V) \subseteq \psi^{-1} (T^2 \times [0,1]) = V$.
Thus we have the following relation:
\begin{equation*}
V = \iota^2 (V) \subseteq \iota (V) \subseteq V ,
\end{equation*}
where we used the fact that $\iota^2$ is trivial.

So we have $\iota (V) = V$,
which means that $\left. \iota \right|_{V}$ is an order two periodic action on $V$
with $\iota ( \partial \E (G) ) = \partial \E (G)$ and $\iota (T) = T$.
Thus we have $\partial (T^2 \times [0,1]) = \psi (\partial \E (G)) \cup \psi(T)$
and $\psi \circ \iota \circ \psi^{-1}$ is an order two periodic action on $T^2 \times [0,1]$
which preserves both $\psi ( \partial \E (G) )$
and $\psi ( T )$.
Then it is shown as Theorem~2.1 in \cite{ms86} that
any finite group action on a three-dimensional manifold
with product structure preserves the structure.
So $\psi ( T )$ must intersect 
with the axis $\psi (V \cap \alpha)$ of the order two periodic action.
This is equivalent to $T \cap \alpha \ne \emptyset$,
which contradicts the assumption that $T \subset \E (G \cup \alpha)$.
Thus $T$ is compressible in $\E (G)$.

We have thus proved Lemma~\ref{lem_parallel}
\end{proof}

\subsection{Simple knots and tunnel number one knots}
A three-dimensional manifold is said to be {\em simple}\/ 
when it contains no essential torus.
A link is said to be {\em simple}\/ 
when it is non-trivial and its exterior is a simple manifold.
By Thurston's study of geometric structures on knot exterior in ${S}^3$,
any non-trivial knot in ${S}^3$ is known to be 
one of a torus knot, a hyperbolic knot or a satellite knot.
See \cite{thp82}; see also \cite{otb01,kab01} for a detailed proof.
Thus any simple knot in ${S}^3$ is either a torus knot or a hyperbolic knot.
For such a knot, we have the following proposition:
%
%
%
%
%
%
\begin{prop} \label{prop_simple}
For any strongly invertible non-trivial simple link $L$ in $S^3$
and the axis $\alpha$ of any strong inversion of the pair $(S^3,L)$,
any embedded equivariant torus in $\E (L)$
which is disjoint from $\alpha$ is compressible in $\E (L)$.
\end{prop}

\begin{proof}
The definition of simpleness for a link $L$ implies that
there is no essential torus in $\E (L)$.
So any embedded equivariant torus in $\E (L)$ is either compressible or boundary-parallel.
However, since $L$ is non-trivial,
the latter case cannot happen by Lemma~\ref{lem_parallel}.
\end{proof}

A knot $K$ in $S^3$ is said to be \textit{tunnel number one}
if there exists an arc $\tau$ properly embedded in $S^3 - K $
such that
$\E ( K \cup \tau )$ becomes a handlebody of genus two.
Such a knot is isotopic onto the genus two surface
appearing as the boundary $\partial \E ( K \cup \tau )$,
which gives a standard genus two
Heegaard surface in $S^3$.
So the knot admits a strong inversion induced from that of the surface.
Thus every tunnel number one knot is known to be strongly invertible.
See also Lemma~5 in \cite{mor95}.
%
%
%
%
%
%
%
\begin{prop} \label{prop_tunnel}
For any non-trivial
tunnel number one knot $K$ in $S^3$
and the axis $\alpha$ of any strong inversion of the pair $(S^3,K)$,
any embedded equivariant torus in $\E (K)$
which is disjoint from $\alpha$ is compressible in $\E (K)$.
\end{prop}
%
%
%
%
%
%
%
%
%
\begin{proof}
Because of Proposition~\ref{prop_simple},
we only consider the case where the tunnel number one knot is satellite.
For such a knot $K$,
it has been proved as Theorem~2.1 in \cite{ms91} that
it is made of
a non-trivial torus knot $K_0$ and a two-bridge link $K_1 \cup K_2$
which is neither a trivial link nor a Hopf link
by gluing $\partial \E (K_2)$ to $\partial \E (K_0)$.
Since $K_2$ itself is a trivial knot in $S^3$,
$\E (K_2)$ is a solid torus containing $K_1$.
Thus the image of $K_1$ 
under this gluing
can be regarded as a knot in $S^3$.
From now on, we will denote the images of $\E(K_1 \cup K_2)$
and $\E (K_2)$ under this gluing
by the same symbols.
Then we say that $\E (K)$ is decomposed into two parts;
one is the torus knot exterior part $\E (K_0)$, which is a Seifert fibered space, and
the other is the two-bridge link exterior part $\E (K_1 \cup K_2)$,
which admits hyperbolic structure. 
Thus the characteristic submanifold of $\E (K)$ coincides with
$\E (K_0)$ up to ambient isotopy in $\E (K)$.
See \cite{jab80}.
Then, by Theorem~8.6 in \cite{ms86},
$\E (K_0)$ can be taken to be preserved by
the strong inversion $\iota$ of the pair $(S^3 , K)$ with axis $\alpha$. 
This implies that the action of $\iota$ gives
a periodic action of order two on $\E (K_0)$,
and it then equals to the restriction to $E(K_0)$ of the strong inversion of the torus knot.
See 10.6 of \cite{ka96b} for example. 

Let $T$ be an embedded equivariant torus in $\E (K)$ which is disjoint from $\alpha$.
If $T$ is not compressible in $\E (K)$,
then one of the following two cases can happen:
\begin{itemize}
	\item The torus $T$ is essential in $\E (K)$, or
	\item The torus $T$ is boundary-parallel in $\E (K)$.
\end{itemize}

If the latter case happens,
then we can apply Lemma~\ref{lem_parallel} so that $K$ must be a trivial knot.
This contradicts the assumption of $K$.
So we study the former case from now on.
Then this assumption implies that $T$ is essential in $\E (K \cup \alpha)$ as well.
Since the satellite torus $T_s$ is incompressible in $\E(K)$, 
the punctured surface $T'_s := T_s \cap \E ( K \cup \alpha )$
is also incompressible in $\E ( K \cup \alpha )$.
Thus, in the former case, $T$ can be isotoped in $\E (K \cup \alpha)$ so that
one of the following two subcases occurs:
\begin{itemize}
	\item[--]
			The torus can be isotoped to $T'$ so that
			the two tori $T'$ and $T_s$ intersect each other,
			and the intersection consists of circles essential in both $T'$ and $T'_s$.
	\item[--]
			The torus can be isotoped to $T'$ so that
			the two tori $T'$ and $T_s$ do not intersect each other.
\end{itemize}

We start with considering the former subcase.
Then we can also assume that the intersection is minimal
up to isotopy in $\E ( K \cup \alpha )$.
Then $T' \cap \E (K_0)$ 
consists of essential annuli in $\E (K_0) - (\E (K_0) \cap \alpha)$.
By taking the projection of $(\E (K_0) , \alpha )$
by the quotient map induced from $\iota$,
we have a two-string tangle $(B^3,t)$. 
\begin{claim}
This tangle $(B^3,t)$ is prime and atoroidal. 
\end{claim}
\begin{proof}
Since $\E (K_0)$ is irreducible and boundary-irreducible, $(B^3,t)$ is prime. 
See Theorem~3.5.15 in \cite{ka96b} for example. 
Since torus knots are simple knots, $\E (K_0)$ is atoroidal.
Thus any embedded equivariant tours in $\E (K_0)$
which is disjoint from $\alpha$ is either boundary-parallel
or compressible in $\E (K_0)$.
If the former case happens,
then $K_0$ must be the trivial knot by Lemma~\ref{lem_parallel}.
So such a torus is always compressible in $\E (K_0)$.
Furthermore,
we can always find such a compression disc in $\E (K_0)$ with no intersection with $\alpha$
by Lemma~\ref{lem_comp}.
This means that
the exterior of $\E (K_0) \cap \alpha$ in $\E (K_0)$ is also atoroidal. 
Then Lemma~\ref{lem_covering} can be applied so that
any embedded torus in $(B^3,t)$ is compressible in $(B^3,t)$,
which means that $(B^3,t)$ is atoroidal.
\end{proof}
On the other hand, 
each essential annulus appearing as a component of $T' \cap \E (K_0)$ 
descends into $B^3 - t$ as an immersed essential, i.e., $\pi_1$-injective, annulus. 
Then we can find an embedded essential annulus in $B^3 - t$
by the Annulus Theorem (see VIII.13 of \cite{jab80} for example).
However this contradicts the following fact: 
if $(B^3,t)$ is prime and atoroidal, then 
$B^3 - t$ contains no essential annulus. 
See 3.6.4 of \cite{ka96b} for example.

Thus the former subcase cannot occur and
now we assume the latter subcase, or
equivalently, 
$T$ is isotoped into $\E (K_0 \cup \alpha)$ or $\E (K_1 \cup K_2 \cup \alpha)$.
Recall that the restriction of the strong inversion $\iota$ into $\E (K_0)$
(resp.\ $\E (K_1 \cup K_2)$) gives its strong inversions with axis $\alpha \cap \E (K_0)$
(resp.\ $\alpha \cap \E (K_1 \cup K_2)$).
Then, since $T$ is essential,
we can apply Proposition~\ref{prop_ett}
so that there is an equivariant essential torus.
This contradicts Proposition~\ref{prop_simple}
since
non-trivial torus knots are known to be simple;
actually they are {\em small}\/ (see Subsection~\ref{subsec_TwoAxes}),
and non-trivial two-bridge links except Hopf link are also known to be simple
(see Corollary~5 of \cite{oe84} for example).

We have thus proved Proposition~\ref{prop_tunnel}.
\end{proof}

These two propositions together with Theorem~\ref{thm_equivalent} \eqref{enu_equatolong}
immediately implies the following result:
%
%
%
%
%
%
\begin{thm} \label{thm_simpletunnel}
Let $K$ be a non-trivial and strongly invertible knot in $S^3$.
If $K$ is simple or tunnel number one,
then, for any $n \geq 3$,
the strongly periodic $\theta_n$-curve
constructed by the $n / 2$-fold cyclic branched covering
along any axis of a strong inversion of the pair $(S^3,K)$ is hyperbolic. \hfill \qed
\end{thm}

As a corollary to Propositions~\ref{prop_simple} and \ref{prop_tunnel},
we have the following result,
which gives a necessary condition for a non-trivial strongly invertible knot to be simple or tunnel number one.

\begin{cor} \label{cor_theta}
Let $K$ be a non-trivial strongly invertible knot in $S^3$
admitting a strong inversion $\iota$ with axis $\alpha$.
Let $G$ be the spatial $\theta_3$-curve
appearing as the quotient of $ K \cup \alpha$ under $\iota$.
Then, if $K$ is simple or tunnel number one,
any embedded torus in $\E (G)$ is compressible.
\end{cor}

\begin{proof}
Take an embedded torus $T$ in $\E (G)$, and
let $\widetilde{T}$ be a torus in $\E( K \cup \alpha )$ appearing as
a component of the pre-image of $T$ by $\iota$.
The torus $\widetilde{T}$ is equivariant under the action of the involution $\iota$,
and disjoint from $\alpha$ in $\E(K)$ by construction.
Since $K$ is assumed to be simple or tunnel number one,
we can find a compression disc for $\widetilde{T}$
which is disjoint from $\alpha$
by Propositions~\ref{prop_simple} and \ref{prop_tunnel} and Theorem~\ref{thm_equivalent}.
Then $\widetilde{T}$ has a compression disc
which is disjoint from $\alpha$ and equivariant under the action of $\iota$
by Lemma~\ref{lem_ECD},
which descends to a compression disc for $T$.
\end{proof}

%
%
%
%
%
%
%
%
%
%
%
\subsection{How about other knots?} \label{subsec_TwoAxes}
Since we have Theorem~\ref{thm_simpletunnel},
remaining strongly invertible knots are those 
which are satellite with tunnel number greater than one.
In this case the property of hyperbolicity of the spatial graphs
generally depends on the choice of the strong inversion.

\begin{figure}[ht]
        \begin{center}
		\includegraphics[width=70mm]{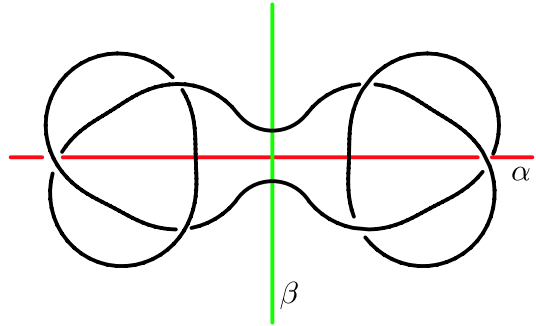}
        \end{center}
        \caption{Connected sum of two trefoil knots with two distinct axes $\alpha$ and $\beta$ of distinct strong inversions}
        \label{f_trefoilsum}
\end{figure}

One of its typical example is a {\em Granny knot}.
This knot is made from a trefoil knot and its copy by their connected sum;
see Figure~\ref{f_trefoilsum}.
A knot in $S^3$ is said to be {\em small}\/
if its exterior contains no closed essential surface.
Since non-trivial torus knots are shown to be small by Theorem in \cite{tsap94},
so is the trefoil knot.
Furthermore, since it is tunnel number one,
we can apply Theorem~4 in \cite{msa00} so that
the Granny knot is not tunnel number one.

The Granny knot $K$ has two distinct strong inversion;
one is with axis $\alpha$,
and the other is with axis $\beta$ in Figure~\ref{f_trefoilsum}.
If we choose the strong inversion with axis $\alpha$,
then we can see that the spatial graphs
obtained by $n/2$-fold cyclic branched covering are hyperbolic,
by checking that the link $\D_{\alpha_1}$,
which appeared in the proof of Theorem~\ref{thm_equivalent},
is hyperbolic.

On the other hand, for $\E (K \cup \beta)$,
one can easily see that 
there is an essential torus,
which is obtained from the decomposing sphere
by tubing along one of the trefoil knot.
This torus is the so-called {\em swallow-follow torus}.
Thus we cannot obtain hyperbolic spatial graphs from this strong inversion
by Theorem~\ref{thm_equivalent}.

Thus we can say that the hyperbolicity depends on the choice of the strong inversion in general.
%
%
%
%
%
%

\section*{Acknowledgements}

The first author is partially supported by Grant-in-Aid for Young Scientists (B),
No.\ 18740038,
Ministry of Education, Culture, Sports, Science and Technology, Japan.

The authors would like to thank Professor Kimihiko Motegi for his encouragement and comments, especially for his comments on Lemma~\ref{lem_comp},
thank Professor Makoto Sakuma for his encouragement and comments, 
especially for Proposition~\ref{prop_tunnel},
and thank Professor Ryo Nikkuni for his comments as a specialist of the spatial graph theory.

The authors also would like to thank the referee for his/her comments and suggestions,
especially for the comment on the proof of Lemma~\ref{lem_comp}.


\providecommand{\bysame}{\leavevmode\hbox to3em{\hrulefill}\thinspace}
\providecommand{\MR}{\relax\ifhmode\unskip\space\fi MR }
\providecommand{\MRhref}[2]{%
  \href{http://www.ams.org/mathscinet-getitem?mr=#1}{#2}
}
\providecommand{\href}[2]{#2}

\end{document}